\documentclass[11pt]{amsart} 
\usepackage{amsfonts,amscd}
\usepackage{psfig}
\newcommand {\partials} {\partial \hspace{-6pt} /}
\newcommand {\partiala} {\partials_{A}}

\newcommand {\stn} {P_{N}}
\newcommand {\qn} {Q_{N}}

\newcommand {\cm} {\chi_{-}}
\newcommand {\al} {\Delta_{N}^{s}}

\newcommand {\cal } {\mathcal}

\def\Z{{\mathbb Z}}
\def\Q{{\mathbb Q}}
\def\R{{\mathbb R}}
\def\C{{\mathbb C}}
\newtheorem{thr}{Theorem}[section]
\newtheorem{df}[thr]{Definition}
\newtheorem{lm}[thr]{Lemma}
\newtheorem{pr}[thr]{Proposition}
\newtheorem{co}[thr]{Corollary} 

\newtheorem{conj}[thr]{Conjecture} 
\begin{document} 
\baselineskip 22pt
\title[Norms on the cohomology of a $3$-manifold and SW theory]{Norms on
the cohomology of a $3$-manifold and SW theory}
\author{Stefano Vidussi}
\address{Department of Mathematics, University of California, Irvine,
California 92697}
\curraddr{Department of Mathematics, Kansas State University, Manhattan, Kansas 66506}
\email{vidussi@math.ksu.edu}
\maketitle
\baselineskip 18pt
\noindent {\bf Abstract.} The aim of this paper is to discuss some
applications of the relation between Seiberg-Witten theory and two
natural norms defined on the first cohomology group of a closed 
$3$-manifold $N$ -
the Alexander and the Thurston norm. We will start by giving a ``new" proof,
applying SW theory, of McMullen's inequality between these two norms, and then
use these norms to study two problems related to symplectic
$4$-manifolds of the form $S^{1} \times N$. First we will prove that - as
long as $N$ is irreducible - the unit balls of the Thurston and
Alexander norms are related in a way that is similar to the case of
fibered $3$-manifolds, supporting the conjecture that $N$ has to be
fibered over $S^{1}$. 
Second, we will provide the first example of a $2$-cohomology class on
a symplectic manifold (of the form $S^{1} \times N$) that lies in the
positive cone and satisfies Taubes' ``more constraints", but cannot be
represented by a symplectic form, disproving a conjecture of Li-Liu 
(\cite{LL2}, Section 4.1).
\section{Introduction}
It has been proven, in \cite{Mc}, that two natural
(semi) norms defined on the first cohomology group of a
$3$-manifold, namely the Alexander
norm  $\|\cdot \|_{A}$, defined from  the Alexander polynomial of the
manifold, and the Thurston  norm $\| \cdot \|_{T}$,  defined in terms of
the minimal genus of the representatives  of the Poincar{\'e} dual two  
dimensional homology  class, satisfy a relation expressed in the following
\begin{thr} \label{tineq} (McMullen) Let $N$ be a compact, connected,
  oriented $3$-manifold (eventually with
  boundary a union of tori); then the Alexander and Thurston norm satisfy 
\begin{equation} \label{ineq} \| \cdot \|_{A}  \leq \| \cdot
\|_{T} + \left\{ \begin{array}{ll} (1 + b_{3}(N)) {\it div}
    ( \cdot) & \mbox{if $b_{1}(N) =
      1$} \\ \\ 0 & \mbox{if $b_{1}(N) > 1$}, \end{array} \right. 
\end{equation} where ${\it div} ( \cdot)$ denotes the divisibility of
    an element in $H^{1}(N,\Z)$. \end{thr} 
This inequality, applied  to  the particular case  where the
three manifold  is the exterior of a  knot $K$,  reduces to the well
known fact that the  degree of  the Alexander  polynomial of the  knot
(i.e.  the difference  between highest and  lowest  power) is bounded from
above by twice  the genus of the  knot,  i.e. the lowest  value of the
genus of a  Seifert  surface of  the knot.  
\\ These two norms turn out to be strictly related 
to the $3$-dimensional Seiberg-Witten theory. 
\\ The proof of Theorem \ref{tineq} given in \cite{Mc} 
is purely topological, but it is suggested the existence of a
proof based on SW theory. (It appears that the first one to observe this has
been D.Kotschick; P.Kronheimer previously proved the inequality in the
case of $N$ obtained as $0$-surgery of a knot, in \cite{K}; the first
detailed proof of the general case appeared in a preprint of
the author (\cite{V}), on which this paper is partly based.)
\\ Our first aim will be to write the two norms in terms of SW basic and
monopole classes for $N$. This allows, as mentioned, an alternative 
proof of Theorem \ref{tineq}, that we will work out for the case, 
for us more interesting, of a closed manifold with $b_{1}(N) >1$.
\\ Then we will use these results to study symplectic $4$-manifolds of the
form $S^{1} \times N$, for an irreducible $N$. We will prove the
following \begin{thr} \label{sameno}
Let $N$ be an irreducible $3$-manifold with $b_{1}(N) > 1$ such that
  $S^{1} \times N$ admits a symplectic structure $\omega$; then 
  there exists a face $F_{T}$ of the unit ball of the Thurston norm
contained in a face $F_{A}$ of the unit ball of the Alexander norm.
\end{thr} 
This quite peculiar property is satisfied by fibered $3$-manifolds, and
supports the conjecture that a $3$-manifold $N$ such that $S^{1} \times N$ admits a symplectic structure must in fact be fibered.
\\ Theorems \ref{tineq} and \ref{sameno} hold true also in the case of
$b_{1}(N) = 1$. It is not surprising that the proof is technically
quite longer, due to the chamber structure of the SW invariants in that
case. We will omit this case here, referring the interest reader to
\cite{V}, where it is treated in detail.
\\ We will then address the problem of determining the constraints for a
cohomology class $\alpha$ on a symplectic $4$-manifold to be represented
by a symplectic form $\omega$. It is clear that such an $\alpha$ must have
positive square, and its pairing with the SW basic classes must satisfy
the constraints determined in \cite{Ta2}. Li-Liu have conjectured in
\cite{LL2} that these are sufficient conditions. We will prove the
following \begin{thr} There exists a symplectic $4$-manifold of the form
$S^{1} \times N$ and a cohomology class $\alpha$ of positive square
satisfying Taubes' ``more constraints" which can not  be represented by a
symplectic form. \end{thr} Hence, Li-Liu conjecture is false.
\section{Alexander and Thurston norms} 
We start by briefly recalling the definition of Alexander
and Thurston norms on the  first cohomology group of a closed,
oriented $3$-manifold  $N$. Denote  by  $F$ the   free abelian  group 
$F  := H_{1}(N,\Z)/Tor$; by definition, $rk(F) = b_{1}(N)$.
The Alexander polynomial of  $N$ is an
element of the group ring  $\Z [F]$, i.e.  a finite sum  
\begin{equation} \label{defal} \Delta_{N} = \sum_{\bf i} a_{\bf i} t^{\bf i}
\end{equation} where
${\bf i} = (i_{1},...,i_{b_{1}(N)})$
is a multi-index of cardinality $b_{1}(N)$, $t =
(t_{1},...,t_{b_{1}(N)})$ with $\{ t_{i} \}$ a basis of $F$ and $a_{\bf
i}$ are integer coefficients. 
The Alexander  polynomial is well defined up  to
multiplication by units of $\Z[F]$. For any element $\phi
\in H^{1}(N;\Z)$ we define the norm \begin{equation} \label{ale} 
\|\phi\|_{A}
:=  \mbox{max}_{\bf ij} \phi(t^{\bf i} \cdot t^{- \bf j}), \end{equation} 
where the indexes run over all $\bf i,j$ such  that $a_{\bf i},a_{\bf j}$
are non  zero. It is clear that this definition is unaffected   by the
indeterminacy in  the Alexander polynomial and does not depend
on the coefficients. 
\\ The Thurston  norm, described in
\cite{Th}, is defined as  follows:  for any Riemann surface
$\Sigma$ embedded in $N$ 
denote  \begin{equation} \cm(\Sigma) = 
\sum_{\Sigma_{i}|g(\Sigma_{i}) \geq 1}
(-\chi(\Sigma_{i})), \end{equation}
where $\Sigma$ is the disjoint union of the $\Sigma_{i}$; 
we then define the norm \begin{equation} \|
\varphi\|_{T}  =     \mbox{min} \{   \cm(\Sigma)   |  \Sigma
\hookrightarrow N, PD[\Sigma]  = \phi \}. \end{equation}
It  is not difficult to verify that
both norms are linear on rays and satisfy the triangle inequality.
It is possible to continuously extend these norms to cohomology
with real coefficients. The unit ball of these norms is then a finite, convex
(possibly noncompact) polyhedron. In particular, the unit ball of the
Alexander norm is by construction dual (up to a factor $2$) to the
Newton polyhedron of $\Delta_{N}$.   
\section{Basic classes and monopole classes} 
In this section we will discuss the way the Alexander and Thurston
norms are related to Seiberg-Witten theory. Essentially
the relation between Alexander norm and SW theory will be deduced from
Meng-Taubes proof of the equivalence of a SW invariant of a $3$-manifold and the Alexander polynomial of the manifold. The relation of Thurston norm and SW theory, instead, has been analyzed in \cite{KM}. 
\\ We start with a brief review of SW theory in dimension three, in
order to have a formulation which is the suitable for our
purposes. Let $(N,g)$ be a
smooth, closed, oriented, riemannian three dimensional manifold. We will
assume that $b_{1}(N) > 1$.  We equip $N$ with the canonical homology
orientation induced by a basis of $F$. 
Once $N$ is endowed with a spin$^{c}$ structure $\stn$, i.e. a
$U(1)$ lifting of the $SO(3)$ frame bundle, we can consider the three
dimensional SW equations \begin{equation} \label{sw3} 
F_{A} = q(\psi) - i \eta, \ \
  \partiala \psi = 0, \end{equation} where $A$ is a connection on the
determinant bundle of the spin$^{c}$ structure, $q(\cdot)$ is an
$\Omega^{2}(N ;i \R)$-valued bilinear form on the sections of the spinor
bundle associated to $\stn$, $\eta$ is a
perturbation term that lives in $\Omega^{2}(N ;\R) \cap ker
\hspace{1pt} d$, and $\partiala$ is the
Dirac operator that acts on spinors. These equations are
invariant under the gauge group of those automorphisms of $\stn$ which act
trivially on the frame bundle. This group acts freely away from
reducible couples, that we can remove suitably choosing good perturbations. 
It is possible to prove (see e.g. \cite{MT}), using standard techniques, that
choosing a generic nonexact perturbation the moduli space of solutions of
Equation \ref{sw3}, modulo gauge equivalence, is a $0$-dimensional compact, 
oriented, smooth manifold; under change of the metric and perturbation 
(as $b_{1}>1$) different moduli spaces are moreover cobordant.
We denote by ${\mathcal M}(\stn,g,\eta)$ the moduli space of solutions of Equations \ref{sw3}, omitting the arguments whenever unnecessary. 
\\ We define the SW invariant for $\stn$ as the algebraic sum of the oriented points of ${\mathcal M}(\stn,g,\eta)$ for $\eta$ a good perturbation.
We have the following definition: \begin{df} Let $c \in
  H^{2}(N; \Z)$ be an integral cohomology class that arises as first
  Chern class of a spin$^{c}$ structure $\stn$ such that the invariant
  $SW(\stn)$ is nonzero. Then $c$ is called a basic class of
  $N$. \end{df} It is quite clear from this definition that the SW
equations for a basic class admit a solution for any metric and a
generic perturbation. Moreover, as the compactness of the equations implies that non emptiness is an open condition, also the unperturbed equations have a solution for any metric, i.e. ${\mathcal M}(\stn,g,0) \neq \emptyset$ (note that this space can be nonsmooth).
This makes it natural to introduce the \begin{df} (Kronheimer-Mrowka) 
Let $c \in H^{2}(N;
  \Z)$ an integral cohomology class that arises as Chern class of a
  spin$^{c}$ structure $\stn$ such that ${\mathcal M}(\stn,g,0) \neq
  \emptyset$ for any metric $g$. Then $c$ is called a monopole class. 
\end{df} From the previous observation, the set of monopole classes,
that we denote by ${\mathcal C}(N)$, contains 
all the basic classes. \\ We now introduce, 
following ref. \cite{MT}, an element in $\Z[[F]]$,
defined from the family of $SW$ invariants of the spin$^{c}$ structures. 
\\ The set ${\mathcal S}$ of spin$^{c}$-structures on $N$ is an affine
$H^{2}(N;\Z)$. There is a natural way to define a map which goes from 
${\mathcal S}$ to $F$ which
is constructed as follows. Fix a reference spin$^{c}$-structure $\qn$, 
that we choose to be the product structure.
Any other structure $\stn$ differs from it by the action of an element
of $H^{2}(N,\Z)$. Consider now the composed map \begin{equation} 
H^{2}(N,\Z) \stackrel{PD}{\longrightarrow} H_{1}(N,\Z)
\stackrel{\pi}{\longrightarrow} F. \end{equation}
Using this map we can construct a map $s$ which goes from ${\mathcal S}$ to $F$.
The fiber of this map is given by the order of the torsion of 
$H_{1}(N; \Z)$, that we will denote now on by $ord(N)$.
Note that twice this map gives, up to torsion, the Poincar{\'e} dual of
the Chern classes of $\stn$. Consider for any $t^{\bf i} \in F$ the
set $s^{-1}(t^{\bf i}) \in {\mathcal S}$. These are the spin$^{c}$
structures that have the same {\it real} Chern class. Define now
\begin{equation} SW(t^{\bf i}) := \sum_{s^{-1}(t^{\bf i})} SW(\stn);
\end{equation} we can now define from this the function \begin{equation}
  \label{ave} SW(N) =
  \sum_{\bf i} SW(t^{\bf i}) t^{\bf i} \in \Z[[F]]. \end{equation} 
Well known facts of SW theory are that the number of
spin$^{c}$ structures for which unperturbed SW equations admit solutions is bounded, and that the invariant $SW(t^{\bf i})$ is symmetric under the natural 
involution of $F$. These observations, together with the definition of the 
function $SW$, yield the fact that $SW(N)$ is a symmetric element of  
$\Z[F]$. 
\par The previous definition, in the case where $ord(N) = 1$, is a simple
reformulation of SW theory. In the other cases, instead, they define
a kind of ``average" over all structures which have the same real
Chern class. 
We can introduce a new definition that is
quite practical for treating the information on spin$^{c}$ structures 
contained in the $SW$ functions of Equation \ref{ave}. 
For any element $\gamma \in H^{2}(N,\Z)$, we denote by
$\gamma^{F}$ its projection to $H^{2}(N,\Z)/Tor (=F^{PD})$.
\begin{df} \label{kmdef} Let $c \in H^{2}(N,\Z)/Tor$ be a cohomology
  class such that \begin{equation} \label{summa} \sum_{c^{F}_{1}(\stn)
      = c} SW(\stn) \neq 0. \end{equation} 
Then $c$ is called an a-basic class (where the ``a" stands for
  averaged). \end{df} We have the following inclusions: 
\begin{equation} \label{incl} 
{\mathcal A}(N) = \mbox{(a-basic classes)} \subset \mbox{(basic classes)}^{F}
\subset \mbox{(monopole classes)}^{F} = {\mathcal C}(N)^{F}. \end{equation}
We want to relate now a-basic classes with the $SW$ function $SW(N)$:
let $c$ be an a-basic class; then the sum appearing in Equation 
\ref{summa} coincides with $SW(t^{\bf i})$ where $t^{\bf i}$ is defined
by  the relation $t^{2 \bf i} = PD (c) $, and the invariant $SW(t^{\bf i})$ 
is nonzero.  
\section{Relation between the norms} 
Our aim now is to relate a-basic classes of $N$ with its Alexander
      polynomial, and then to the Alexander norm.
In this section we will give a proof of the following 
\begin{pr} \label{prop} 
Let $N$ be a closed three manifold with $b_{1}(N) > 1$; then
      the Alexander norm of an element $\phi \in H^{1}(N;\Z)$ is given
      by \begin{equation} \| \phi \|_{A}  =  \mbox{max}_{{\mathcal A}(N)}
(c \cdot \phi) \end{equation}  where the maximum is taken over all
a-basic classes of $N$. \end{pr}
\noindent {\bf Proof}: The basic ingredient for the proof is provided by the
      theorem of Meng and Taubes which identifies the $SW$ function
      with the (sign-refined) Reidemeister-Franz torsion introduced by
      Milnor. This is related, on its own, to the (sign-refined)
      symmetrized Alexander polynomial, denoted by $\al$.
More precisely, we have the \begin{lm} (Meng-Taubes, Turaev)
      Let $N$ be a closed three manifold equipped with its canonical
      homology orientation, with $b_{1}(N) > 1$; then we
      have, in $\Z[F]$,
 \begin{equation} \label{swm} SW(N) = \al. \end{equation}
       \end{lm}
From this relation, considering the definitions in Equation \ref{defal} and \ref{ave}, we have $SW(t^{\bf i}) = a_{\bf i}$ and, by Definition
\ref{kmdef}, a-basic classes $c_{\bf i} \in {\cal A}(N)$ are twice
the Poincar\'e duals of elements of $F$ with nonvanishing invariant $SW(t^{\bf i})$, i.e. \begin{equation} \label{obeq} c_{\bf i} \in H^{2}(N; \Z)/Tor 
\mbox{ is a-basic} \Longleftrightarrow a_{\bf i} \neq 0 \ \ \mbox{where} 
\ \ t^{2 \bf i} = PD(c_{\bf i}). 
\end{equation} The use of the relations of Equation
\ref{swm} 
allows us to write the Alexander norm in terms of a-basic classes: we can
write \begin{equation} \label{dualpo}
 \| \phi \|_{A} = \mbox{max}_{\bf ij} \phi(t^{\bf i} \cdot t^{- \bf j})
= \mbox{max}_{{\mathcal A}(N)} (c \cdot \phi). \ \
  \end{equation} 
The latter equality comes as follows: for any fixed couple
${\bf i,j}$ with nonzero coefficient we have
\begin{equation} \label{refe} 
\phi(t^{\bf i} \cdot t^{- \bf j}) = \phi(t^{\bf i}) -
\phi(t^{j})
\leq  \mbox{max}_{\bf ij}(|\phi(t^{2 \bf i})|,|\phi(t^{2 \bf j})|) \leq
\mbox{max}_{\bf k}|\phi(t^{2 \bf k})|
\end{equation} where ${\bf k}$ ranges among all indexes with nonzero 
coefficient. Being $\Delta_{N}^{s}$ symmetric, 
equality in Equation \ref{refe} is attained for some choice of index with 
${\bf j} = - {\bf i} = \pm {\bf k}$. 
By Equation \ref{obeq}, we can thus write
$\phi(t^{2 \bf k}) =  (c_{\bf k} \cdot \phi)$ (and remove the absolute value).
This completes the proof of our statement. \endproof
The content of Proposition \ref{prop} is the good one for our
purpose, because of the results of \cite{KM} on the relations between monopole
classes and Thurston norm: we have now all we need to prove Theorem
\ref{tineq}.
\\ First, we observe that we can restrict the proof to the case
of an irreducible manifold. In fact, the SW polynomial vanishes
for the connected sum of manifolds with $b_{1} > 0$, and the 
set of a-basic classes of a manifold is preserved for connected sum with a
rational homology sphere (the $SW$ polynomial gets multiplied by the
order of the torsion of the rational homology sphere).
\\ For irreducible manifolds we have now the equality expressed in   
\begin{thr} \label{kromro}
  (Kronheimer-Mrowka) Let $N$ be a manifold as above: then the
  Thurston norm of a class $\phi \in H^{1}(N,\Z)$ is given by 
\begin{equation} \label{krom} \| \phi \|_{T} = \mbox{max}_{{\mathcal
C}(N)}  (c \cdot \phi). \end{equation} \end{thr} 
Putting together the inclusion ${\mathcal A}(N) \subset {\mathcal
C}(N)^{F}$, Theorem \ref{kromro} and Proposition \ref{prop}, we deduce
the inequality \ref{ineq} for any closed three manifold
(for the sole purpose of proving inequality \ref{tineq}, it is sufficient
at this point to prove an inequality in Equation \ref{krom} for all basic 
classes, as in \cite{A}).
\\ {\bf Remark:} Equation \ref{dualpo} states that the unit ball of the 
Alexander norm, up to a factor,  is dual to the Newton polyhedron of $SW(N)$. 
Equation \ref{krom}, in light of the content of \cite{KM} (see also \cite{K}) 
states that the Thurston unit ball is dual, up to the same factor, 
to the polyhedron of ``$SWF(N)$'', 
an element of $\Z[F]$ that can be constructed from the Seiberg-Witten-Floer 
invariants of $N$ (still lacking a rigorous treatment). This answers a 
longstanding question of Fried in \cite{Fr}.
\section{Symplectic $S^{1} \times N$} \label{sympls}
In this section we will use the Alexander and Thurston norm to study the
following conjecture: \begin{conj} \label{conj} (Taubes) Let $N$ be a $3$-manifold
such that
$S^{1} \times N$ admits a symplectic structure $\omega$. Then $N$ admits
a fibration over $S^{1}$. \end{conj}
We will assume again that $b_{1}(N) > 1$.  Under this condition
manifolds that fiber over
$S^{1}$ are irreducible, and it is known that for any class $\phi \in
H^{1}(N,\Z)$ representing a fibration we have $\| \phi  \|_{A} = \| \phi
\|_{T}$ (see e.g. \cite{Mc}). 
Fibered classes are known to satisfy the following condition 
(see \cite{Th}): the integral points laying in the cone over a 
(top dimensional) face of the Thurston unit sphere have the property of 
being all fibered, or none does.
This implies in particular that a fibered face of the unit
ball of Thurston norm is contained in a face of the unit ball of the
Alexander norm. \\ We would like to prove that the latter condition holds
for any $N$ such that $S^{1} \times N$ is symplectic.
We will be able to do so under the further assumption that
$N$ is irreducible; in view of the results of \cite{McC}, this is a
reasonable assumption. Our proof adapts to the case of $b_{1} > 1$ 
the strategy of \cite{K1}.
\\ We observe that, as simplecticity is an open condition, there is no
restriction  (see
\cite{D}) in assuming that the symplectic form $\omega$ on $S^{1} \times
N$ is the reduction of an integer class of $H^{2}(S^{1} \times N,\Z)$. 
There is a cone, in $H^{2}(S^{1} \times N,\R)$, of cohomology classes that
can be represented by symplectic forms, and in this cone the set of classes
which are in the image of the cohomology with rational coefficients is
dense. We will be interested to have cohomology
classes that lie in the image of the cohomology with integral
coefficients, and eventually pass to sufficiently high multiples of the
symplectic form: we will implicitly assume this whenever necessary.
\\ We want to recall now some general results that we will 
apply to our case. The first is the Donaldson theorem on the existence
of symplectic submanifolds (\cite{D}). 
This theorem assures that there exist a connected symplectic 
submanifold $H \subset S^{1}
\times N$ such that \begin{equation} \label{subma} 
[H] = PD[\omega] = [S^{1}] \times \gamma +
\tau \in H_{2}(S^{1}
\times N, \Z) \end{equation} where $\gamma \in H_{1}(N,\Z)$,
$\tau \in H_{2}(N,\Z)$ and $\gamma \cdot
\tau > 0$ (for sake of
notation we will denote all products, both on $N$ and on $S^{1} \times N$,
with a dot, the distinction being clear from the context). Denote $\phi =
PD(\tau) \in H^{1}(N,\Z)$; as a consequence of the previous discussion,
the $\phi$'s associated to symplectic forms as in the relation of Equation
\ref{subma} define a cone in $H^{1}(N,\Z)$.
\\ The second result is that the spin$^{c}$
structures on $S^{1} \times N$ with nontrivial $SW$ invariants 
must be pull backs of spin$^{c}$ structures on $N$ (to prove this you can
use, e.g., the adjunction inequality); moreover there is an
identification between the moduli spaces  for a spin$^{c}$ structure
$\stn$  on $N$ and the moduli space for the
pull-back structure on $S^{1} \times N$ (that we will usually denote
with the same symbol), once a suitable correspondence of the
perturbation terms is set (see \cite{OT}). This allows the identification, up
to a sign determined by the choice of homology orientations, of the SW
invariants associated to these moduli spaces.
\\ The third point concerns spin$^{c}$ structures on a symplectic four
manifold $(M,\omega)$ with canonical bundle $K$. There exist, in
that case, a canonical spin$^{c}$ structure that decomposes as $\C
\oplus K^{-1}$ (and has first Chern class equal to $-K$). 
Any other spin$^{c}$ structure can be written as $E \oplus (K^{-1}
\otimes E)$ for an $E \in H^{2}(M,\Z)$. 
There are some constraints on spin$^{c}$ structures with nonvanishing
invariants that arises from Taubes' work (see \cite{Ta1},\cite{Ta2}).
In the case of $b_{+}(M) > 1$ 
the canonical spin$^{c}$ structure has $SW$ invariant $\pm 1$ and
for any other structure $E_{i} \oplus
(K^{-1} \otimes E_{i})$ with nonzero invariants we have 
$ K \cdot \omega  \geq E_{i} \cdot \omega \geq 0$. Equality implies, 
respectively, $E_{i} = \C$ or $E_{i} = K$. \\ This inequality translates,
for the basic classes $\kappa_{i} := det(E_{i} \oplus
(K^{-1} \otimes E_{i}))$, in the relation \begin{equation} \label{taco} K
\cdot \omega \geq | \kappa_{i} \cdot \omega |, \end{equation} with
equality only for the case
$\kappa_{i} = \pm K$.
Let's apply these results to $S^{1}
\times N$. First, the canonical class and all other basic
classes are pull backs: there exists
a preferred line bundle $K \in H^{2}(N,\Z)$ (for sake of simplicity, we
use the same notation on
$N$) and a preferred spin$^{c}$ structure on $N$ of the form $\C \oplus
K^{-1}$ with $SW$ invariant $\pm 1$ such that any other spin$^{c}$
structure on $N$ appears as $E_{i} \oplus (K^{-1} \otimes E_{i})$ for 
$E_{i} \in H^{2}(N,\Z)$. The structures with nonzero invariants  must satisfy 
$K \cdot \phi \geq E_{i} \cdot \phi \geq 0$, the equalities implying
respectively $E_{i} = \C$ or $E_{i} = K$. 
This translates to a constraint, for the basic classes
of $N$, which has the form  \begin{equation} \label{threeco} K
\cdot \phi \geq | \kappa_{i} \cdot \phi |,
\end{equation} with equality only for the case
$\kappa_{i} = \pm K$.
\\ Using this it is straightforward to prove the following \begin{pr} 
\label{prev} Let $(S^{1} \times N,\omega)$ be a symplectic manifold with
$b_{1}(N) > 1$, and denote by $\phi \in H^{1}(N,\Z)$ the K\"unneth
component of $[\omega]$; then $\|\phi\|_{A} = K \cdot \phi$. 
Moreover, $\phi$ lies in the cone over a top dimensional face of
the unit ball of the Alexander norm, dual to the vertex $K$ of the
Newton polyhedron of $\Delta_{N}$. \end{pr}
\noindent {\bf Proof}: The maximum of 
$\kappa_{i} \cdot \phi$ for $\kappa_{i}$ basic is attained for and
only for $K$. We want to use this property to evaluate the
Alexander norm, in the form expressed in Proposition
\ref{prop}. To do this we need only to prove that $K$ (or, more
precisely, its image $K^{F}$ in $H^{2}(N,\Z)/Tor$) is an a-basic class. 
But no other basic class $\kappa_{i}$ can coincide up to torsion with
$K$ without violating Equation \ref{threeco}, so that the sum of Equation
\ref{summa}, namely $\sum_{c^{F}_{1}(\stn) =  K^{F}} SW(\stn)$, contains
only one nonzero term, that term being equal to $1$. 
This means that $K^{F}$ is an a-basic class. We can conclude,
following  Proposition \ref{prop}, that 
$\|\phi\|_{A} = K \cdot \phi$. The rest of the proposition is an
obvious consequence of what was previously stated.
\endproof
We will use Proposition \ref{prev} to write the genus 
of the symplectic submanifold $H$ of
equation \ref{subma}, in conjunction with the adjunction inequalities
for manifolds of type $S^{1} \times N$ that are contained in \cite{K1}. 
These
apply to irreducible manifolds $N$ which do not have a basis of $H_{2}(N,\Z)$
composed of tori. Leaving aside this totally degenerate case, 
for which the equality of Alexander and Thurston norm is trivial, we have
the following 
\begin{pr} \label{iden} Let $(S^{1} \times N,\omega)$ be a symplectic 
manifold with $N$ irreducible,
$b_{1}(N) > 1$, and denote by $\phi \in H^{1}(N,\Z)$ the K\"unneth
component of $[\omega]$: then $\|\phi\|_{A} = \|\phi\|_{T}$. \end{pr} 
{\bf Proof}: The adjunction inequality for embedded submanifolds of
$S^{1} \times N$ of \cite{K1} can be written in the form \begin{equation}
  \chi_{-}(H) \geq H \cdot H + \| \phi \|_{T} =  2 \gamma \cdot \tau +
\| \phi \|_{T}. \end{equation} As $H$ is symplectic, the adjunction
formula for symplectic submanifolds gives 
\begin{equation} \label{chi} \chi_{-}(H) = H \cdot
  H  + K \cdot H = 2 \gamma \cdot \tau + \| \phi \|_{A} \end{equation}
 These formulae are compatible
with the content of Equation \ref{tineq}
if and only if $\|\phi\|_{A} = \|\phi\|_{T}$. \endproof
We can somehow strengthen this result. By Equation \ref{tineq} the unit ball of
the Thurston norm is contained in the unit ball of the Alexander norm; it
is clear that extending the Alexander norm to real coefficients, and
using the denseness of $H^{1}(N,\Q)$, the equality stated in Proposition
\ref{iden} continues to hold in an open cone of $H^{1}(N,\R)$ determined
by the cone of classes of $H^{2}(S^{1} \times N,\R)$ admitting a
symplectic representative (the norm is a continuous function). Therefore
a (top dimensional) face $F_{T}$ of the unit ball of the Thurston norm 
(containing $\phi/\|\phi\|_{T}$) intersects a face $F_{A}$ of the unit ball 
of the Alexander norm (the face dual to $K$); but this implies that
the entire $F_{T}$ is contained in $F_{A}$. 
\begin{figure}[h] 
\centerline{\psfig{figure=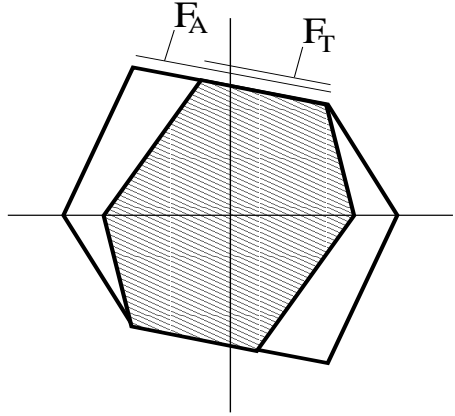,height=55mm,width=60mm,angle=0}}
\caption{\label{sam} {\sl Alexander unit ball and Thurston unit ball for a $3$-manifold such that $S^{1} \times N$ is symplectic.}} 
\end{figure}
Figure \ref{sam} shows a possible 
case of the relation between the norms for $N$ as described.
\\ This observation completes
the proof of Theorem \ref{sameno}.
\section{Representability of a cohomology class by a symplectic form} 
A classical problem of symplectic topology is to determine necessary and
sufficient conditions for a $2$-cohomology class
$\alpha$ on an even-dimensional closed, smooth, oriented manifold $M$ 
admitting an almost complex structure with canonical class $K$, to be 
represented by a symplectic form. 
This problem merges with the
general problem of the existence of {\it any} symplectic structure
on $M$. Necessary conditions arise
from the very definition of symplectic form; in particular, if $M$ has
dimension $4$, we need $\alpha \in {\mathcal P}$, where
\begin{equation} {\mathcal P} = \{ \beta \in H^{2}(M,\R) | \beta \cdot
\beta > 0 \}. \end{equation} An early conjecture, in \cite{Thsy}, speculated 
that every almost complex manifold with nonempty ${\mathcal P}$ admitted 
a symplectic structure. Since then, other constraints have been identified.
In particular, more refined conditions arise from
Taubes' constraints on SW basic classes: as mentioned, we must have 
$SW(K) = \pm 1$, and the class $\alpha$ must satisfy the conditions of
eq. \ref{taco}.  We denote by ${\cal T}$
the cone composed of elements of ${\cal P}$ satisfying these constraints,
i.e. \begin{equation} \label{tauco} 
{\mathcal T} := \{ \alpha \in {\mathcal P} | K \cdot
\alpha \geq |\kappa_{i} \cdot \alpha | \}, \end{equation}
with strict inequality when $\kappa_{i} \neq \pm K$. 
\\ It is well known that satisfying Taubes' constraints is not a
sufficient condition for $\alpha$ to be represented by a symplectic form.
In fact, as discussed in
\cite{KMT}, if we  consider the manifold $X \# \Sigma$ where $X$ is
symplectic and $\Sigma$ is an homology $4$-sphere admitting a nontrivial
cover, and the cohomology class $\alpha_{\omega}$ on $X \# \Sigma$
induced by a symplectic form $\omega$ on $X$ under the natural
isomorphism 
$H^{2}(X,\R) = H^{2}(X \# \Sigma,\R)$, we have identity of the
Seiberg-Witten polynomials $SW_{X} = SW_{X \# \Sigma}$ and
$\alpha_{\omega}$ lies in ${\cal T}_{X \# \Sigma}$,
but cannot be represented by a symplectic form for the simple
reason that $X \# \Sigma$ itself does not admit symplectic structures
(it has a cover with trivial SW polynomial).  There is another class of
potential, more refined, examples of couples $(M,\alpha)$ satisfying these
constraints with $\alpha$ not representable by a symplectic form. These
are knot surgery manifolds homotopic to a $K3$ surface (see \cite{FS} for
the definition) obtained from a knot $K$: it is commonly conjectured that
whenever $K$ is not a fibered knot, $M$ can not be symplectic,
but it is easy to find nonfibered knots such that 
Taubes' constraints are satisfied for a class $\alpha$. 
\\ In both the previous cases, the absence of a symplectic
form representing $\alpha$ has to be attributed, in some sense, to the
manifold $M$ (which does not admit {\it tout court} symplectic
structures) and not to the cohomology class itself. We can ask about the
situation for manifolds known to be symplectic. In particular, it has
been conjectured (see 
\cite{LL2}, Section 4) that if we assume that $X$ is a symplectic
manifold, the cone ${\mathcal T}$
 coincides with the
``symplectic cone" \begin{equation} {\mathcal W} := \{ \alpha \in
{\mathcal P} | \alpha \mbox{ is represented by a symplectic form }
\}. \end{equation} The conjecture gives a possible answer to the
following problem, outlined in the beginning of this section, namely
\\ {\bf Question}: Let $M$ be a symplectic manifold. Determine the cone,
in $H^{2}(M,\R)$, represented by symplectic forms. 
\\ Some partial answer to this question are known. For example, Geiges
proved in \cite{Ge} that for $T^{2}$-bundles over $T^{2}$, all classes in the positive
cone are represented by symplectic forms (we remark that all these
classes satisfy Taubes' constraints, as the canonical class is trivial);
it is interesting also to compare with the result of Gromov for the case of
open manifolds, where {\it any} form in the positive cone lies in
${\mathcal W}$.
\\ Concerning this Question, we have the following result: \begin{thr}
\label{dem} There exist symplectic manifolds, of the form $S^{1} \times
N^{3}$,  on which there are cohomology classes of positive
square satisfying Taubes' constrains but which can not be 
represented by symplectic
forms, i.e  the strict inclusion ${\mathcal W} \subset {\mathcal T}$
holds true. In particular, the conjecture of \cite{LL2} is
false. \end{thr} {\bf Proof}: The proof is based on the following
assumption, that will be proved in the next section (Theorem \ref{lastt}): 
There exists a
family of fibered $3$-manifolds, whose generic component is denoted by
$N$, with $H_{1}(N,\Z) = \Z^{2}$, such that the fibered face
$F_{T}$ of the Thurston unit ball is strictly contained in the face
$F_{A}$ of the Alexander unit ball. Assuming this, we proceed as follows. 
Denote by ${\cal V}$
the (nonempty) cone, in $H^{1}(N,\R)$, over $F_{A} \setminus {\bar
F_{T}}$. Choose a $\phi \in {\cal V}$. We claim that we can define an a
class $\psi \in H^{2}(N,\R)$ such that the cohomology class $\alpha
\in H^{2}(S^{1} \times N, \R)$ with K\"unneth decomposition
\begin{equation} \alpha = \phi \wedge [dt] + \psi \end{equation}
has positive square and satisfies Taubes' constraints, i.e. $\alpha
\in {\cal T}$. This is achieved in the following way. 
We have $\alpha
\cdot \alpha = 2 \phi \cdot \psi$: identify $H^{2}(N,\R) =
Hom(H^{1}(N,\R),\R)$; to get a positive square, we can choose $\psi$ to be 
any element of the cone $Hom(\phi,\R_{+})$. 
We observed in Section \ref{sympls} that the basic classes
on $S^{1} \times N$ are pull-back of basic classes on $N$; the choice of
$\psi$ is therefore irrelevant for the constraints of Equation \ref{tauco} 
and $\alpha$ belongs to ${\cal T}$ if and only if $\phi$ satisfies 
the condition (on $H^{*}(N,\R)$) $K
\cdot \phi \geq | \kappa_{i} \cdot \phi |$
with equality only for the case $\kappa_{i} = \pm K$. But this condition
is equivalent to the condition that $\phi$ lies in the cone over $F_{A}$, as 
$F_{A}$ is, by definition, the face dual to $K$, i.e. the elements lying in the cone over $F_{A}$ have maximal pairing, among all basic classes, 
with and only with $K$. 
\\ To complete the proof, we need to show now that $\alpha$ can not be
represented by a symplectic form. By Proposition \ref{iden} (and the
following comments, if we want to work with 
cohomology with real coefficients), if
$\alpha$ admits a symplectic representative, then its K\"unneth component
$\phi$ should have the same Alexander and Thurston norm, something we
excluded choosing $\phi \in {\cal V}$. \endproof Note that proceeding as
above we can without difficulty choose the class $\alpha$ to lie in the
image of cohomology with integer coefficients.
\\ {\bf Remarks:} 1. The symplectic manifold discussed in Theorem \ref{dem} 
is not simply
connected, but we believe that there exist simply connected examples. 
In particular, we expect that the link surgery manifolds obtained using a link with fibered face strictly contained in a face of the Alexander norm 
(as the one we will discuss in the next section), are possible examples. 
The difficulty in proving such result arises from the difficulty of proving 
the 
analogue of Proposition \ref{iden} (see Section 7 of \cite{K} for a 
discussion of this interesting problem). 
\\ 2. The failure of Li-Liu conjecture, as expressed in the examples of 
Theorem \ref{dem}, is due to the mismatch between the convex hull of basic 
classes and the convex hull of monopole classes, as the latter determines the 
extension of the fibered cone of $H^{1}(N,\R)$. It is coincevable to improve 
the conjecture, at least for symplectic manifolds of the form 
$S^{1} \times N$, by reformulating the definition of the cone 
${\mathcal T}$ as \begin{equation} {\mathcal T} := \{ \alpha \in {\mathcal P} | K \cdot \alpha \geq |\kappa_{i} \cdot \alpha | \mbox{\ \ for any  
monopole class $\kappa_{i}$} \}, \end{equation}
with strict inequality when $\kappa_{i} \neq \pm K$. In that case, as follows 
from the results of \cite{Th}, the 
conjecture would hold true assuming the validity of a strict version of 
Conjecture \ref{conj}, which takes the form \begin{conj} 
Let $N$ be a $3$-manifold such that $S^{1} \times N$ admits a symplectic 
structure $\omega$; then the K\"unneth component of $[\omega]$ in 
$H^{1}(N,\R)$ can be represented by a nondegenerate $1$-form (i.e. it lies in 
a fibered cone). \end{conj}
\section{Construction of the three manifolds} 
\begin{figure}[h] 
\centerline{\psfig{figure=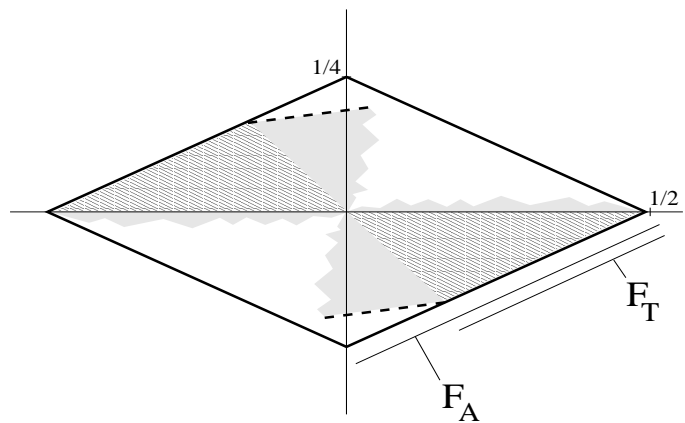,height=55mm,width=90mm,angle=0}}
\caption{\label{mynorm} {\sl Alexander unit ball and Thurston unit ball for 
$N(K_{1},K_{2})$ with $g(K_{1}) = 2, g(K_{2}) = 4$ - the dotted regions are qualitative.}} 
\end{figure}

In this section we will justify the assumption made in the proof of
Theorem \ref{dem}, namely the existence a family of closed, fibered
$3$-manifolds with the property that $F_{T}$ is strictly contained in
$F_{A}$. Our construction  will be based on the existence of a
noteworthy $2$-component link, exhibited by Dunfield in \cite{Du}, which
has the same property. We will need the following result: \begin{pr}
\label{dunf} (Dunfield)
There exists a $2$-component oriented link $D = D_{1} \cup D_{2} \subset
S^{3}$ with Alexander polynomial \begin{equation} \Delta_{D}(t_{1},t_{2})
= (t_{1} - 1)(t_{2} - 1) \end{equation} (written in terms
of the homology classes of the meridians to the two components) which has
a fibered face $F_{T}$ strictly contained in a face $F_{A}$, dual to the 
vertex $t_{1}t_{2}^{-1}$ of the dual polyhedron. \end{pr} (In
Dunfield paper, the Alexander polynomial and the norms are discussed in
terms of an homology basis different from ours, but it is easy to rewrite
them in terms of the standard homology basis for the link exterior, as
above). \\ We don't know the exact shape of the Thurston unit ball, but
for our purpose it is enough to  know the result contained in the
previous Proposition. Denote $N_{D} = S^{3} \setminus \nu (D_{1} \cup
D_{2})$:  
\begin{figure}[h] 
\centerline{\psfig{figure=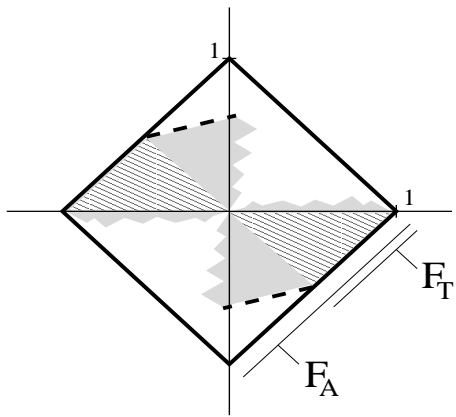,height=55mm,width=60mm,angle=0}}
\caption{\label{figdunf} {\sl Alexander unit ball and Thurston unit ball for $N_{D}$ - the dotted regions are qualitative.}} 
\end{figure}
Figure
\ref{figdunf} represents, in the space $H^{1}(N_{D},\R)$ with
basis vectors the dual basis $\tau_{i}$, the unit ball
of the Alexander norm, and a part of the unit ball of the Thurston norm.
\\ Let now $K_{1}$, $K_{2}$ be a
couple of fibered knots of genus $g(K_{i}) > 0$ and let
$\Delta_{K_{i}}(t)$ be their Alexander polynomials. Next, define the closed
manifold \begin{equation} N(K_{1},K_{2}) = N_{D} \cup (\coprod_{i=1}^{2}
S^{3} \setminus \nu K_{i}), \end{equation} where on the boundary tori the
gluing map is defined to be the orientation reversing diffeomorphism
which identifies the basis $(\mu(D_{i}),\lambda(D_{i}))$ with
$(\mu(K_{i}),-\lambda(K_{i}))$. To interpret this, notice that each knot
exterior is an homology solid torus, so that this operation appears
as an homology Dehn filling for $N_{D}$, with surgery
coefficient $0$. The reason of the choice of this surgery curve appears
evident from the fact (for a proof, see \cite{EN}, Section 3) that the
minimal genus Seifert surface in $N_{D}$ representing an homology class
Poincar\'e  dual to a class $(m_{1},m_{2}) \in H^{1}(N_{D},\Z)$
intersects the boundary torus $T_{i}$ in  $m_{i}$ copies of the longitude
of each link component (note that $lk(D_{1},D_{2}) = \Delta_{D}(1,1) = 0$). 
Each minimal genus Seifert surface of $N_{D}$ has
therefore a natural capping in $N(K_{1},K_{2})$, given by the union
of $m_{i}$ copies of the fiber of $S^{3} \setminus \nu K_{i}$. 
In particular, if $(m_{1},m_{2})$ is a fibered class, this fibration
extends to a fibration of $N(K_{1},K_{2})$ through the fibrations
\begin{equation} S^{3} \setminus \nu K_{i} \longrightarrow S^{1}
\stackrel{( \cdot )^{m_{i}}}{\longrightarrow} S^{1} \end{equation} of the
knots' exteriors. This proves, in particular, that $N(K_{1},K_{2})$ is
irreducible. As the linking number of $D$ is zero,
$H_{1}((N(K_{1},K_{2}),\Z) = \Z^{2}$, canonically identified with
$H_{1}(N_{D},\Z)$.
\\ It is natural to guess that, for any class $(m_{1},m_{2}) \in
H^{1}(N_{D},\Z)$, the surface constructed above is the minimal genus
representative for the cohomology class Poincar\'e dual to
$(m_{1},m_{2}) \in H^{1}(N(K_{1},K_{2}),\Z)$.  
This is confirmed by the following \begin{lm}
\label{eisn} (Eisenbud-Neumann, Prop. 3.5): If $M$ is a compact
irreducible
$3$-manifold and ${\bf m} \in H^{1}(M)$, then the Thurston norm $\|{\bf
m}\|_{T}$ is the sum of the norms of the restrictions of $\bf m$ to the
Jaco-Shalen-Johannson components of $M$. \end{lm} As a consequence of
this  Lemma, denoting with the symbol $\| \cdot \|_{\hat T}$ the norm on
the closed manifold, we have the following
\begin{co}  The Thurston norm $\|(m_{1},m_{2})\|_{\hat T}$ of an element 
$(m_{1},m_{2}) \in H^{1}(N(K_{1},K_{2}),\Z)$ is given by \begin{equation}
\label{tuno}
\|(m_{1},m_{2})\|_{\hat T} = \|(m_{1},m_{2})\|_{T} +  |m_{1} | (2 g
(K_{1}) - 1) + |m_{2} | (2 g(K_{2}) - 1) \end{equation} where
$\|(m_{1},m_{2}) \|_{T}$ is the Thurston norm of the corresponding
element of $N_{D}$. \end{co} {\bf Proof}: this follows from Lemma
\ref{eisn}, together with the observation that the class $(m_{1},m_{2})$
on the closed manifold restricts to the element with same coordinates in
$N_{D}$ and to the classes $m_{i} \in H^{1}(S^{3} \setminus \nu
K_{i},\Z)$, which have Thurston norm $\|m_{i}\|_{T} = |m_{i}|(2g(K_{i})
-1)$, by definition of genus of a knot and linearity on rays. \endproof 
We want to study now the Alexander norm of the manifold $N(K_{1},K_{2})$;
in order to do this we need a gluing formula for the Alexander
polynomial (or the SW invariant) along tori. We have the
following \begin{lm} (Gluing formula) Let $N(K_{1},K_{2}) = N_{D} \cup
(\coprod_{i=1}^{2} S^{3} \setminus \nu K_{i})$ be defined as above: then 
the Alexander polynomials of the manifolds are related by the formula
\begin{equation} \label{poly} \Delta_{N(K_{1},K_{2})}(t_{1},t_{2}) =
\Delta_{N_{D}}(t_{1},t_{2}) \frac{\Delta_{K_{1}}(t_{1})}{t_{1} - 1}
\frac{\Delta_{K_{2}}(t_{2})}{t_{2} - 1} = \Delta_{K_{1}}(t_{1})
\Delta_{K_{2}}(t_{2}). \end{equation} Therefore, the Alexander norm $\|
\cdot \|_{\hat A}$ on $N(K_{1},K_{2})$ is given by \begin{equation}
\label{alno} \begin{array}{c}
\|(m_{1},m_{2})\|_{\hat A}  = |m_{1} | 2g(K_{1}) +
|m_{2} | 2g(K_{2})  = \\ \\ = \|(m_{1},m_{2})\|_{A} +  |m_{1} | (2 g
(K_{1}) - 1) + |m_{2} | (2 g(K_{2}) - 1).  \end{array} \end{equation}
\end{lm} {\bf Proof}: It is known (see \cite{MT}, \cite{T2})
that the Milnor torsion is multiplicative by gluing along tori, with
suitable identification of the variables; this torsion coincides with the
Alexander polynomial for manifolds having $b_{1} > 1$ and it is equal to
the Alexander polynomial $\Delta_{K}(t)$ divided by $(t-1)$ for the case
of a knot.  
Remembering Proposition \ref{dunf}, Equation
\ref{poly} above follows. The relation on the Alexander norm is then an
easy corollary of this formula, as the degree of the Alexander
polynomial of a fibered knot equals twice its genus.
\endproof This Lemma says, in particular, that the unit ball of the
Alexander norm for $N_{D}$ and $N(K_{1},K_{2})$ are conformally
equivalent (see Figure \ref{mynorm}). We are ready to prove the
\begin{thr} \label{lastt} 
There exist a family of fibered closed $3$-manifolds $N$ with
$H_{1}(N,\Z) = \Z^{2}$ such that a fibered face 
of the Thurston unit ball is strictly contained in the corresponding face
of the Alexander unit ball. \end{thr} {\bf Proof:} Our family is
given, for any choice of the fibered knots $K_{i}$, by $N(K_{1},K_{2})$. 
We observed before that each fibration of $N_{D}$ in $F_{T}$ extends to a
fibration of $N(K_{1},K_{2})$, defining a fibered face $F_{\hat T}$ of the
Thurston unit ball for $N(K_{1},K_{2})$. This face will be contained in
$F_{\hat A}$, one of the four faces of the Alexander unit ball 
(having the same cone as $F_{A}$). This face is dual to the vertex 
$t_{1}^{2g(K_{1})}t_{2}^{-2g(K_{2})}$ (square of a vertex of the Newton 
polyhedron of the symmetrized Alexander polynomial).
If a class $(m_{1},m_{2}) \in
H^{1}(N_{D},\Z)$ lies in the cone over $F_{A} \setminus {\bar F_{T}}$ (in
particular
$\|m_{1},m_{2}\|_{A} < \|m_{1},m_{2}\|_{T}$), then the corresponding
class on the closed manifold has Alexander norm strictly smaller than
the Thurston norm, from Equations \ref{tuno} and \ref{alno}, i.e. $F_{\hat T}$
is strictly contained in $F_{\hat A}$. From this the statement follows.
\endproof 
\noindent Figure \ref{mynorm} describes the Thurston and Alexander norm for a 
particular choice of $g(K_{i})$.
\\ We want to outline a second proof of the same
statement, based on the fact that a class on a closed three manifold is
fibered if and only if all the restrictions to each JSJ component are
fibered (see
\cite{EN}, Theorem 4.2): If the Thurston norm on the closed manifold
coincided with the Alexander norm on a larger cone than the one on
$N_{D}$, then there would be fibered classes on $N(K_{1},K_{2})$ which 
restrict, on $N_{D}$, to
nonfibered ones (as mentioned above, 
all integral points laying on a face of the unit ball of
the Thurston norm containing at least one fibration are fibered).
\\ We finish this section pointing out that, although the link $D$ above is 
the only example worked out in detail, fibered links with the properties 
of Proposition \ref{dunf} are likely to be ``frequent'' (compare the discussion
in \cite{Du}). From these examples, other closed $3$-manifolds can be 
constructed.

\end{document}